\documentclass[11pt]{amsart}

\usepackage{amsmath}
\usepackage{amssymb}
\usepackage{graphicx}
\usepackage{amscd}
\usepackage{stmaryrd}
\usepackage[mathscr]{euscript}
\theoremstyle{definition}

\numberwithin{equation}{section}

\newcommand{\sB}{\mathscr B}
\newcommand{\sC}{\mathscr C}

\newcommand{\sH}{\mathscr H}
\newcommand{\sJ}{\mathscr J}

\newcommand{\sM}{\mathscr M}
\newcommand{\sN}{\mathscr N}
\newcommand{\sP}{\mathscr P}
\newcommand{\sR}{\mathscr R}

\title[Condenser Quasicentral Modulus]{The condenser quasicentral modulus  }

\author[D.-V. Voiculescu]{Dan-Virgil Voiculescu${}$}

\address{Department of Mathematics \\ University of California at Berkeley \\ Berkeley, CA\ \ 94720-3840}
\email{{\tt dvv@math.berkeley.edu}}

\begin{document}

\begin{abstract}
We introduced the quasicentral modulus to study normed ideal perturbations of operators. It is a limit of condenser quasicentral moduli in view of a recently noticed analogy with capacity in nonlinear potential theory. We prove here some basic properties of the condenser quasicentral modulus and compute a simple example. We also discuss some associated noncommutative variational problems. Part of the results are in the more general setting of a semifinite von Neumann algebra.
\end{abstract}

\maketitle

%%%%%%%%%%%%%%%%%%%%%%%%%%%%%%%%%%%%%%%%%%%%%%%%%%%%%%%%%%%%%%%%%%%%%%%%%%%%%%%%

\section{Introduction}
\label{sec1}
The quasicentral modulus (\cite{8},\cite{9},\cite{10}) plays a key role in the study of Hilbert space operators modulo normed ideals (see our surveys \cite{11},\cite{12}).This paper is a sequel to \cite{13}.In \cite{13} we made the case that the quasicentral modulus is a noncommutative analogue of capacity in nonlinear potential theory, where the first order Sobolev spaces use general rearrangement invariant norms of the gradients. One consequence is that the quasicentral modulus becomes a limit quantity of condenser quasicentral moduli. Note that
the condenser quasicentral moduli are usually finite and non-zero also in situations where the quasicentral modulus can take only the values $0$ and $ \infty $, like in the case of the $p$ - classes when $p > 1$. Another new feature is that we will often deal with the more general case of separable semifinite factors or von Neumann algebras, that is, not only with the type I case of the algebra  $\sB (\sH)$ of bounded operators on the Hilbert space $\sH$ and note also that even the case when $\sH$ is finite dimensional is no longer a trivial case.

Concerning the nonlinear potential theory capacity with which we observed an analogy see \cite {2} and more references are in \cite{13}. Prior to this, we had noticed connections with Yamasaki hyperbolicity (\cite{14}) and with the noncommutative potential theory based on Dirichlet forms \cite{1}.

\bigskip
\noindent
\underline{\qquad\qquad} 

2020 {\em Mathematics Subject Classification}. Primary: 46L89; Secondary: 31C45, 47L20.

{\em Key words and phrases}. quasicentral modulus, noncommutative nonlinear condenser capacity, normed ideals of operators, semifinite von Neumann algebras, noncommutative variational problems.

\vfill
\noindent

\newpage
Can the analogy with nonlinear potential theory be further extended ? A way to achieve this may be via noncommutative variational problems . We take a small first step in this direction introducing variational problems related to the condenser. We observe that computations in the case of the p-classes naturally lead to noncommutative analogues of the p-Laplace equation.

Besides the introduction and references there are eight more sections. Section~2 is about preliminaries and basic definitions. Section~3 contains some general properties of the condenser quasicentral modulus in the semifinite setting. In particular we prove that under certain conditions the the condenser quasicentral modulus is symmetric with respect to switching the projections which define the condenser. We also give a result about the behavior with respect to certain conditional expectations. Section ~4 gives a lower bound for the condenser quasicentral modulus in the case of the algebra of bounded operators on a Hilbert space. This is analogous to the lower bound in \cite{10} for the quasicentral modulus. Section ~5 is the computation of an example arising from the bilateral shift operator. In section ~6 we adapt and generalize to our semifinite setting the result in \cite{9} about the largest reducing projection on which the quasicentral modulus vanishes. In the analogy with nonlinear potential theory, this is a special noncommutative polar set.Section ~7 deals with variants of the quasicentral condenser modulus. Section ~8 is about noncommutative variational problems. We make some general remarks about minimizers  for the condenser problem.

\section{Preliminaries and Definitions}
\label{sec2}

We introduce here the framework in which we will work, especially related to normed ideal / symmetric operator space norms (\cite{3}, \cite{4}, \cite{5}, \cite{6}) 
and we recall the definition of the quantities we introduced in \cite{13}.

By $(\sM, \rho)$ we will denote a von Neumann algebra $\sM \subset \sB(\sH)$ where $\sH$ is a separable complex Hilbert space and $\rho$ is a faithful normal semifinite trace on $\sM$. We will assume that $\sM$ is either atomic, that is it is generated by its minimal projections or that it is diffuse and $\rho(I) = \infty$, in particular there are no minimal projections in this case. Thus $\sM$ could be for instance a type ${I}$ or a type ${II}_\infty$ factor with its trace, but it could also be for instance $L^\infty (S, d\sigma)$ where the measure sigma has no atoms and is not finite or it could be $\ell^\infty (X)$ with the measure giving mass $1$ to each singleton subset ( the traces are those corresponding to the measures).

We will denote by $Proj(\sM)$ the selfadjoint projections in $\sM$ and by $\sP(\sM)$ or simply $\sP$ the set of $P \in Proj(\sM)$ so that $\rho(P) < \infty$. By $\sR$ we will denote the set of $x \in \sM$ for which there is $P \in \sP$ so that $xP = x$ , that is the ideal of operators of finite $\rho -rank $ in $\sM$. Further, $\sR^+_1$ will stand for the positive contractions in $\sR$ , that is $\{a\in \sR \vert 0 \le a \le {I} \}$.
It will also be convenient to introduce the set $\Lambda = \{ L \in L^1 (\sM , \rho) \vert 0 \le L \le {I} \}$. Thus we have
$\sP \subset \sR^+_1 \subset \Lambda \subset \sM$.

If $x \in \sM$, the generalized singular values ( see \cite{5}) are
\begin{center}
\[
\mu(t, x) = \inf \{ \Vert A({I} - P) \Vert \vert P \in \sP, \rho (P) le t \}
\]
\end{center}
where $t>0$  and $\mu (x) $ will denote the function 
\begin{center}
\[
(0, \infty)  \longrightarrow \mu(t,x) \in \lbrack 0 , \infty ).
\]
\end{center}

On $L^1(\sM,\rho) \cap \sM$ we will consider a norm $\vert \cdot \vert_\sJ$ to which we will refer as the symmetric operator norm. We will assume that $\vert \cdot \vert_\sJ$ satisfies the folllowing conditions:
\begin{itemize}
\item [1.] 
 $\mu(x) \le \mu(y) \Rightarrow \vert x \vert_\sJ \le \vert y \vert_\sJ$
which has among its consequences, that if $a, b \in \sM, x\in \sM \cap L^1(\sM, \rho)$ then
$\vert axb \vert_\sJ \le \Vert a \Vert \vert x\vert_\sJ \Vert b \Vert $.
\item [2.]
$ C_1 \min (\vert x \vert_1 , \Vert x \Vert ) \le \Vert x \Vert _\sJ \le C_2 ( \vert x \vert_1 + \Vert x \Vert )$
for some constants $C_1, C_2 \in (0 , \infty)$.
\item [3.]
If $P_n \in \sP$ are so that $\rho(P_n) \rightarrow 0$ as $n \rightarrow \infty$ then $\vert P_n \vert_\sJ \rightarrow 0$. 
This condition is trivially satisfied if $\sM$ is a type ${I}$ factor, because it is meaningless there being no non-zero $P_n$ as above. The condition is equivalent to $x_n \in \sM\cap L^1 (\sM , \rho), \Vert x_n \Vert \le 1, \vert x_n\vert_1 \rightarrow 0 \Rightarrow \vert x_n \vert_\sJ \rightarrow 0 $.  
The condition can also be put in the form : there is an increasing function $\phi : (0, \infty) \longrightarrow (0, \infty)$ so that 
$\lim_{t\rightarrow0}  \phi(t) = 0$ for which we have $\Vert x\Vert \le 1, x \in L^1(\sM , \rho) \Rightarrow \vert x \vert_\sJ \le \phi(\vert x\vert_1 )$.
Note also that if $(\sM , \rho)$ is diffuse and the condition is not satisfied, then there must be a constant $C \in (0, \infty)$ so that $\vert x \vert_\sJ   \ge C \Vert x\Vert $.
\item [4.]
Considering the von Neumann algebra $\sM \otimes \frak M_n$ endowed with the trace $\rho \otimes Tr_n$, the norm $\vert \cdot \vert_\sJ$ on $\sM \cap L^1(\sM, \rho)$, identified with a subspace of $ \sM \otimes e_{11} $ has an extension to a norm, we will still denote by  $\vert \cdot \vert_\sJ $ on  $(L^1 (\sM , \rho) \cap \sM) \otimes \frak M_n) = (L^1 (\sM \otimes \frak M_n , \rho \otimes Tr_n) \cap (\sM \otimes \frak M_n)$ satisfying the analogues of 1. - 3. . 
\end{itemize}

We opted for this ad hoc way of introducing $\vert \cdot \vert_\sJ$, instead of a discussion starting with operator spaces \cite{5} which would have taken us farther than the more modest aim of this paper. We will also denote by $\sJ$ the completion of $L^1 (\sM ,\rho) \cap \sM $ with respect to the norm $\vert \cdot \vert_\sJ$, which is an $\sM $ bimodule. In the case of $\sB(\sH)$, $\sJ$ identifies with a normed ideal in $\sB (\sH) $ consistent with the notation we used in previous papers (see \cite{12} ).  Note also that property 4. in the symmetric operator spaces setting , actually follows from the general relation of symmetric function spaces, symmetric sequence spaces and symmetric operator spaces (see 2.5 in \cite{5} , Thm.2.5.3 and the discussion of Questions 2.5.4 and 2.5.5 in \cite{5} )

\bigskip
\noindent
{\bf Definition 2.1}
Let $P, Q \in \sP, PQ = 0 $ and let $\tau = (T_j)_{1 \le j \le n}, T_j \in \sM, 1 \le j \le n$. The condenser quasicentral modulus with respect to the symmetric operator norm $\vert \cdot \vert_\sJ $ is the number:
\[
{k_\sJ (\tau ; P, Q) = \inf \{ \max_{1 \le j \le n}   | [T_j , A] |_\sJ  \vert  A\in \sR^+_1 , AP = P, AQ = 0 \} } .
\]
Similarly if $\alpha = (\alpha_j)_{1 \le j \le n}$ is a $n$-tuple of automorphisms of $\sM$ which preserve $\rho$ we define
\[
{k_\sJ (\alpha ; P, Q) = \inf \{ \max_{1 \le j \le n}  | \alpha_j (A) - A |_\sJ \vert A \in \sR^+_1 , AP = P, AQ = 0 \} } .
\]

\bigskip
We remark that if $P_1, Q_1 \in \sP, P_1 Q_1 = 0$ and $P \le P_1, Q \le Q_1$ then
\[
{ \{A \in \sR^+_1 \vert AP = P, AQ = 0\}  \supset \{A \in \sR^+_1 \vert AP_1 = P_1 , AQ_1 = 0 \}}
\]
so that 
\[
{k_\sJ (\tau ; P, Q) \le k_\sJ ( \tau; P_1, Q_1)} ,
\]
\[
{k_\sJ ( \alpha ; P, Q) \le k_\sJ (\alpha ; P_1, Q_1)} .
\]

\bigskip
\noindent
{\bf Definition 2.2}
Let $P, Q \in Proj (\sM), PQ = 0$ and let $ \tau, \alpha, \vert \cdot \vert_\sJ $ be like in Definition 2.1 . We define
\[
{k_\sJ (\tau; P, Q) = \sup \{ k_\sJ (\tau ; P', Q' ) \vert P', Q' \in \sP, P' \le P, Q' \le Q \} }
\]
\[
 {k_\sJ (\alpha ; P , Q) = \sup \{ k_\sJ (\alpha ; \vert P', Q' \in \sP, P' \le P, Q' \le Q \} }.
 \]
 Further we define
 \[
 {k_\sJ(\tau ; P ) = k_\sJ (\tau ; P , 0 )}
 \]
 \[
 {k_\sJ (\alpha ; P ) = k_\sJ ( \alpha ; P , 0) }
\]
and the quasicentral moduli of $\tau$ and of $ \alpha$ 
\[
{ k_\sJ ( \tau) = k_\sJ (\tau ; I ) }
\]
\[
{ k_\sJ (\alpha) = k_\sJ ( \alpha : I ) } .
\]
If $\sM$ needs to be specified we will write $ k_{\sJ , \sM} ( \tau ; P, Q ) $ and so on.

\bigskip
In case $(\sM , \rho) = (\sB (\sH) , Tr )$ , this definition of $k_\sJ (\tau)$ is equivalent to the definition we used in our earlier work,
as we already pointed out in \cite{13}.

It will also be useful to have a technical result about replacing $\sR^+_1$ by the larger set $\Lambda$ in the above definitions.
 
\bigskip
\noindent
{\bf Lemma 2.1}
Let $P, Q \in \sP, \tau, \alpha$ be as in Definition 2.1.. We have
\[
{ k_\sJ (\tau ; P, Q) = \inf \{ max_{1 \le j \le n} |[T_j ; A] | _\sJ \vert A \in \Lambda, AP = P, AQ = 0 \} } .
\]
\[
{k_\sJ (\alpha ; P ,Q) = \inf \{max_{1 \le j \le n} | \alpha_j (A) - A |_\sJ | A \in \Lambda , AP = P, AQ = 0 \} }.  
\]

\bigskip
\noindent
{\bf {\em Proof. } }
The Lemma is an immediate consequence of the fact that 
\[
{\Xi = \{ A \in \sR^+_1 \vert AP = P, AQ = 0 \} }
\]
is a dense subset with respect to the topology defined by  the norm $ \vert \cdot \vert_\sJ $ of the set
\[
{\Theta = \{ A \in \Lambda \vert AP = P, AQ = 0 \} . }
\]
This in turn is seen as follows:
\[
{ \Xi = P + ( I - P - Q ) \sR^+_1 (I - P - Q) }
\]
\[
{ \Theta = P + ( I - P - Q ) \Lambda  ( I - P - Q ) }
\]
so that the proof reduces to the proof of the density of $\sR^+_1$ in $ \Lambda$ . If  $ X \in \Lambda $ then let $P_j$ be the spectral projection of $X$ for $ [1/j , 1 ] $.
We then have $ P_j \in \sP $ and $ \Vert XP_j - X \Vert  \rightarrow 0 , | XP_j - X |_1 \rightarrow 0 $ so that $XP_j \in \sR^+_1 $ and by condition 2. we have $ | X - XP_j|_\sJ \rightarrow 0$.
\qed

\section {Some general properties }
\label{sec3}  
\bigskip
\noindent
{\bf Proposition 3.1}
Assume $k_\sJ(\tau) = 0$ and $P, Q \in Proj (\sM)$ are so that $PQ = 0$. Then
\[
{k_\sJ(\tau; P,Q) = k_\sJ(\tau; Q, P) .}
\]
Similarly if instead of $\tau$ we have $\alpha$ with $ k_\sJ (\alpha) = 0 $, then
\[
{k_\sJ(\alpha; P, Q) = k_\sJ(\alpha; Q, P). }
\]

\bigskip
\noindent
{\bf {\em Proof.} }
Since $ k_\sJ (\cdot ;P, Q)$ for general $P, Q$ is the $sup$ of such quantities with $P, Q \in \sP$, it suffices to prove the Proposition under the additional assumption that $P, Q \in \sP$. Moreover, by symmetry it clearly suffices to prove that
\[
{k_\sJ (\cdot ; P, Q) \ge k_\sJ (\cdot ; Q, P)}
\].
If $\epsilon > 0$ there is 
\[
{A \in P+ (I - P - Q) \sR^+_1 (I - P - Q)}
\]
so that 
\[
{ \max_{1 \le j \le n} | [ A , T_j] |_\sJ \le k_\sJ (\tau ; P , Q) + \epsilon .}
\]
Then we have 
\[
{\max_{1 \le j \le n} |[(I - A) , T_j]|_\sJ = \max_{1 \le j \le n} |[A , T_j]|_\sJ \le k_\sJ (\tau ; P , Q) + \epsilon . }
\]
Clearly $0 \le I - A \le I , (I - A)P = 0, (I - A )Q = Q$ and we would be done if there weren't the problem that $I - A$ is not in $\sR$ in general. Since $k_\sJ (\tau) = 0$ we also have $k_\sJ (\tau ; P+Q) = 0$, so that there is $B \in \sR^+_1$ so that  $B(P+Q) = P+Q$ and $\max_{1 \le j \le n} |[B , T_j|_\sJ < \epsilon . $ Let $F = B(I - A)B $. We have:
\[
{F \in \sR^+_1 , FP = 0, FQ = Q}
\]
and
\[
{ |[F , T_j]|_\sJ \le 2 |[B , T_j]|_\sJ + |[I - A , T_j]|_\sJ < 3 \epsilon . }
\]
Thus 
\[
{k_\sJ (\tau ; P,Q) + 3 \epsilon \ge k_\sJ (\tau ;Q, P) }
\]
and $\epsilon > 0$ being arbitrary, $3 \epsilon $ is as good as $\epsilon$ here.

The case of the $n$-tuple of automorphisms $\alpha$ is dealt along the same lines.
If $A \in P + (I - P - Q) \sR^+_1 (I - P - Q) $ and 
\[
{ \max_{1 \le j \le n} |\alpha_j (A) - A |_\sJ \le k_\sJ (\alpha ; P , Q) + \epsilon}
\]
then we have
\[
{ |\alpha_j (I - A) - (I - A)|_\sJ = |\alpha_j (A) - A|_\sJ}
\]
and $ (I - A)P = 0 , (I - A)Q = Q . $
Choosing $B$ as in the previous case, we consider $ F = B (I - A) B$ and we have
\[
{ |\alpha_j (F) - F|_\sJ \le 2 |\alpha_j (B) - B|_\sJ + |\alpha_j (A) - A|_\sJ .} 
\]
This leads then to 
\[
{k_\sJ (\alpha ; P, Q) + 3 \epsilon \ge k_\sJ (\alpha ; Q, P) }
\]
and so on. 
\qed
\bigskip

The framework for the next result involves a von Neumann subalgebra $\sN \subset \sM$ so that $\rho | \sN$ is semifinite, in which case $\sR (\sN) $ is weakly dense in $\sN$. Let $E$ be the conditional expectation of $\sM$ onto $\sN$ with respect to $\rho$ ( see \cite{7} Prop 2.36). If the $n$-tuple of automorphisms $\alpha$ is so that $\alpha_j(\sN) = \sN , 1 \le j \le n $ , then in view of the assumption $\rho \circ \alpha_j = \rho, 1 \le j \le n$ we will have
that $\alpha_j|\sN \circ E = E \circ \alpha_j , 1 \le j \le n $. We shall also assume that 
\[
{x \in L^1 (\sM) \cap \sM \Rightarrow |E(x)|_\sJ \le |x|_\sJ .}
\]
Since the symmetric operator norm $| \cdot |_\sJ $ was introduced in an ad-hoc way, without going into the theory of spaces of operators, we prefer to treat this as an assumption.

\bigskip
\noindent
{\bf Proposition 3.2}
Let $\sN$ be a von Neumann subalgebra of $\sM$ so that $\rho | \sN$ is semifinite and let $E$ be the conditional expectation of $\sM$ onto $\sN$ with respect to $\rho$. Assume also that the $n$-tuple of automorphisms $\alpha$, which preserve $\rho$ is so that $\alpha_j (\sN) = \sN, 1 \le j \le n $ . Let further $P, Q \in Proj (\sN) $ be so that $PQ = 0$. Then we have
\[
{k_{\sJ,\sM} (\alpha ; P,Q) = k_{\sJ,\sN} (\alpha | \sN ; P, Q) .} 
\]
\bigskip
\noindent
{\bf {\em Proof.}}
We first prove the statement in case $P, Q \in \sP (\sN)$. The inequality 
\[
{k_{\sJ, \sN} (\alpha | \sN ; P,Q) \ge k_{\sJ, \sM} (\alpha ; P, Q) }
\]
is obvious, because the $LHS$ is an $inf$ over a subset of the set the $inf$ of which is the $RHS$,

We have $E (\Lambda (\sM)) = \Lambda (\sN) $, so that:
\[
{E (P + (I - P - Q) \Lambda(\sM) (I - P - Q) ) = P + (I - P - Q) \Lambda(\sN) (I - P - Q). }
\]
Thus, if 
\[
{A \in P + ( I - P - Q) \Lambda (\sM) (I - P - Q)}
\]
then we have
\[
{|\alpha_j(E(A) - E(A)|_\sJ = |E(\alpha_j (A) - A) |_\sJ \le |\alpha_j (A) - A|_\sJ }
\]
and hence
\[
{k_{\sJ , \sN} (\alpha | \sN ; P, Q) \le k_{\sJ , \sM} (\alpha ; P , Q)} .
\]
This concludes the proof in the case of $P, Q \in \sP (\sN)$.

Assume now that we only have $P, Q \in Proj (\sN)$. Then we get  that 
\[
\begin{aligned}
{k_{\sJ, \sM} (\alpha;P , Q) = sup \{ k_{\sJ, \sM} (\alpha; P_1 ;Q_1) | P_1 \le P, Q_1 \le Q, P_1, Q_1 \in \sP (\sM) \}} \\
{\ge sup \{ k_{\sJ , \sM} (\alpha ; P_1, Q_1) | P_1 \le P, Q_1 \le Q , P_1, Q_1 \in \sP (\sN) \}} \\
{= sup \{ k_{\sJ , \sN} (\alpha | \sN ; P_1, Q_1) | P_1 \le P, Q_1 \le Q, P_1 , Q_1 \in \sP (\sN) \}} \\
{= k_{\sJ, \sN} (\alpha | \sN ; P, Q) .}
\end{aligned}
\]

We still must prove that if $P_1, Q_1 \in \sP(\sM), P_1 \le P, Q_1 \le Q$ then there are $P_2 , Q_2 \in \sP (\sN), 
P_2 \le P, Q_2 \le Q $ so that 
\[
{ k_{\sJ , \sM} (\alpha ; P_1, Q_1) \le k_{\sJ , \sM} (\alpha ; P_2, Q_2) + \eta }
\]
for a given $ \eta > 0 $.

If $A \in P_2 + (I - P_2 - Q_2) \Lambda (\sM) (I - P_2 - Q_2) $ where $P_2 \le P, Q_2 \le Q, P_2, Q_2 \in \sP(\sN)$ we consider
\[
{B = P_1 + (I - P_1 - Q_1) A (I - P_1 - Q_1)}
\]
so that 
\[
{B \in P_1 + ( I - P_1 - Q_1) \Lambda (\sM) ( I - P_1 -Q_1) . }
\]
Since $\rho | \sN$ is semifinite, there are $ P_2, Q_2 \in \sP (\sN), P_2 \le P, Q_2 \le Q$ so that 
\[
{ | ( I - P_2 ) P_1|_1 < \epsilon , |(I - Q_2) Q_1|_1 < \epsilon  }
\]
which also implies
\[
{| P_1 (I - P_2) |_1 < \epsilon, |Q_1 (I - Q_2) |_1 < \epsilon }.
\]
Then we have
\[
\begin{aligned}
{ | (I - P_1 - Q_1)(I - P_2 - Q_2) - (I - P_2 - Q_2 ) |_1} \\
{= | (P_1 + Q_1)(I - P_2 - Q_2)|_1  = | P_1 (I - P_2) + Q_1 (| - Q_2)|_1 < 2 \epsilon .}
\end{aligned}
\]
If $ X \in \sM , \| X \| \le 1$ this gives
\[
{| ( I - P_1 -Q_1) (I - P_2 -Q_2) X (I - P_2 - Q_2)(I - P_1 - Q_1) - (I - P_2 - Q_2) X (I - P_2 - Q_2)|_1 < 4 \epsilon }
\]
so that in particular
\[
\begin{aligned}
{ | B - A|_1 \le 4 \epsilon + | P_1 - P_2 + ( I - P_1 - Q_1) P_2 (I -P_1 - Q_1) |_1} \\
{= 4 \epsilon + |P_1 - P_2 + (I - P_1) P_2 (I - P_1)|_1 = 4 \epsilon + |P_1 - P_1 P_2 - P_2 P_1 + P_1 P_2 P_1 |_1} \\
{\le 5 \epsilon + | P_1 P_2 P_1 - P_2 P_1|_1 = 5 \epsilon + | - P_1 (I - P_2) P_1 + (I - P_2 P_1)P_1 |_1 \le 7 \epsilon . }
\end{aligned}
\]

We shall now use the fact that $ |\alpha_j (x) |_\sJ = |x|_\sJ $ which is a consequence of $\mu (\alpha_j (x)) = \mu (x) $
which in turn follows from $ \rho \circ \alpha_j = \rho $.

On the other hand 
\[
{|\alpha_j (B) - B|_\sJ - |\alpha_j (A) - A|_\sJ \le |\alpha_j ( A - B) - (A - B)|_\sJ \le 2 |A - B|_\sJ .}
\]
In view of condition \  3. satisfied by $ | \cdot |_\sJ $ we have
\[
{| A - B |_\sJ \le 2 \phi ( 2^{-1} | A - B|_1 ) \le 2 \phi (4 \epsilon) . }
\]
This in turn gives 
\[
{ k_{\sJ , \sM} (\alpha ; P_1 , Q_1 ) \le k_{\sJ , \sM } (\alpha ; P_2 , Q_2 ) + 4 \phi (4 \epsilon)}
\]
which concludes the proof .
\qed 

\bigskip
\noindent
There is also a similar result for $n$-tuples $\tau$ of self-adjoint operators instead of the $n$-tuple $\alpha$ of automorphisms. Since the proof is along the same lines, we will leave out many details.

\bigskip
\noindent
{\bf Proposition 3.3}
Let $\sN$ be a von Neumann subalgebra of $\sM$ so that $\rho |\sN$ is semifinite and let $E$ be the conditional expectation of $\sM$ onto $\sN$ with respect to $\rho$. Let further $\tau$ be a $n$-tuple of selfadjoint elements in $\sN$ and let $P, Q \in  Proj (\sN) $ be so that $PQ = 0$. Then we have
\[
{k_{\sJ , \sM} (\tau; P, Q) = k_{\sJ, \sN} (\tau; P, Q). }
\]

\bigskip
\noindent
{\bf {\em Sketch of Proof.}}
We first deal with $P, Q \in \sP (\sN)$ . Obviously we have
\[
{k_{\sJ , \sN} (\tau ; P, Q) \ge k_{\sJ , \sM} (\tau ; P, Q). } 
\]
On the other hand $E(\Lambda(\sM)) = \Lambda(\sN) $ and if 
\[ 
{A \in P+ (I - P - Q) \Lambda (\sM) (I - P - Q)}
\]
then
\[
{ | [ T_j , E(A)] |_\sJ = |E([T_j , A])|_\sJ \le |[T_j , A]|_\sJ }
\]
which gives
\[
{ k_{\sJ , \sN} (\tau ; P, Q) \le k_{\sJ , \sM} (\tau ; P,Q) .}
\]
This concludes the proof in case $P, Q \in \sP (\sN) $.

If we only have $P, Q \in Proj (\sN)$ then the preceding immediately gives
\[
{k_{\sJ, \sM} (\tau; P, Q) \ge k_{\sJ, \sN} (\tau; P, Q) .}
\]
We still must show that if  $P_1, Q_1 \in \sP (\sM) , P_1 \le P , Q_1 \le Q$ then there are $P_2, Q_2 \in \sP(\sN) , P_2 \le P, Q_2 \le Q $ so that 
\[
{k_{\sJ, \sM} (\tau ; P_1, Q_1) \le k_{\sJ, \sN} (\tau ;P_2, Q_2) + \eta }
\]
for a given $\eta > 0 $.

Since $\rho | \sN$ is semifinite, there are $P_2, Q_2 \in \sP (\sN) $ so that $P_2 \le P, Q_2 \le Q$ and
\[
{|(I - P_2) P_1|_1 < \epsilon , |(I - Q_2) Q_1|_1 < \epsilon . }
\]

If 
\[
{A \in P_2 + ( I - P_2 - Q_2) \Lambda (\sM) (I - P_2 - Q_2) }
\]
we consider
\[
{ B = P_1 + (I - P_1 - Q_1) A (I - P_1 - Q_1) \in P_1 + (I - P_1 - Q_1 ) \Lambda (\sM) (I - P_1 - Q_1) . }
\]
Then we have 
\[
{ | B - A |_1  \le 7 \epsilon  }
\]
so that
\[
{ |B - A|_\sJ  \le 2 \phi ( 2^{-1} |B - A |_1) \le 2 \phi (4 \epsilon) .}
\]
It follows that
\[
{|[\tau , B] |_\sJ - |[\tau, A]|_\sJ \le 2 \Vert \tau \Vert \phi (4 \epsilon) .}
\]
This in turn shows that 
\[
{k_{\sJ ,\sM} (\tau ; P_1, Q_1) \le k_{\sJ , \sN} (\tau; P_2, Q_2) + 2 \Vert \tau \Vert \phi (4\epsilon) . }
\]
An appropriate choice of $\epsilon > 0 $ concludes the proof.
\qed

\bigskip
\noindent
Taking $P = I $ and $Q = 0$ in Prop. 3.2. and Prop 3.3. we have the following.

\bigskip
\noindent
{\bf Corollary 3.1}
Let $\sN$ be a von Neumann subalgebra of $\sM$ so that $\rho | \sN$ is semifinite. If $\alpha = (\alpha_j)_{1 \le j \le n} $ is a $n$-tuple of automorphisms which preserve $\rho$ , $\alpha_j (\sN) = \sN , 1 \le j \le n$ and $\tau = (T_j)_{1 \le j \le n} $ is a $n$-tuple of selfadjoint elements of $\sN$, then we have
\[
{ k_{\sJ, \sM} (\alpha) = k_{\sJ, \sN} (\alpha | \sN) ,  k_{\sJ, \sM} (\tau) = k_{\sJ, \sN} (\tau) .}
\]

\bigskip
\noindent
{\bf Remark 3.1} 
If $\sM = \sM_1 \otimes \sM_2 $ , where $(\sM_1 , \sigma)$ is a type ${II_1}$-factor and $\sM_2$ is a factor of type ${I_\infty}$ so that
$\rho = \sigma \otimes Tr $ and if $\tau = I \otimes \tau_2$ where $\tau_2$ is a $n$-tuple of selfadjoint operators in $\sM_2 = \sB (\sH_2)$ for some Hilbert space $\sH_2$, the preceding corollary gives that
\[
{k_\sJ (\tau) = k_{\sJ, \sM_2} (\tau_2) . }
\]
This shows in particular that the examples of quasicentral modulus of $n$-tuples in $\sB (\sH_2)$ give automatically  examples of quasicentral modulus in type ${II_\infty}$ factors by taking $I \otimes \tau_2 $ in $\sM_1 \otimes \sB(\sH_2) $. One can proceed in a similar way for $n$-tuples of automorphisms ( for the type ${I_\infty}$ case, the next proposition shows that this reduces to the quasicentral modulus for $n$-tuples of unitary operators) .

\bigskip
\noindent
Let us also make a very simple observation about the case of $n$-tuples of unitary operators. If $ u = (U_j)_{1 \le j \le n} $ is a 
$n$-tuple of unitary elements of $\sM$ we denote by  $ Ad u = (Ad U_j)_{1 \le j \le n} $ the $n$-tuple of inner automorphisms where 
$(Ad U_j) (x) = U_j x U^*_j $. Consider also a map $\epsilon : \{ 1, ... , n\} \longrightarrow \{ 1 , \ast \} $ and let then $ u^{\epsilon} =       
( U^{\epsilon (j)}_j )_{1 \le j \le n} $. \\ We have 
\[
\begin{aligned}
{ |[U_j , A] |_\sJ = |[U_j , A] U^{\ast}_j |_\sJ = | (Ad U_j ) (A) - A|_\sJ} \\ {=|U^{\ast}_j [U_j , A] |_\sJ = |(Ad U^{\ast}_j ) (A) - A|_\sJ = |[U^{\ast}_j , A ] |_\sJ .}
\end{aligned}
\]
This immediately implies the following Proposition.

\bigskip
\noindent
 {\bf Proposition 3.4}
 Let $ u = (U_j)_{1 \le j \le n} $ be a $n$-tuple of unitary elements of $\sM$ , let $\epsilon: \{1, ... , n\} \longrightarrow \{ 1, \ast \}$ and let $P, Q \in Proj (\sM ), PQ = 0$. Then we have
 \[
 {k_\sJ (u; P, Q) = k_\sJ (u^{\epsilon}; P, Q) = k_\sJ ( Ad u^{\epsilon} ; P, Q). } 
\]
 
\section{The lower bound}                                                                                                                                                                                                                        
\label{sec4}
In this section we assume $\sM = \sB(\sH)$ and that $|\cdot |_\sJ$ is the symmetric norm arising from a norming function $\Phi$, so that we will write $|\cdot |_{\Phi}, k_{\Phi} (...) $ instead of $|\cdot |_\sJ , k_\sJ (...) $. By $\Phi^{\ast} $ we will denote the dual norming function, so that if $\frak S^{(0)}_{\Phi}$ is the closure of $\sR$ in the norm $|\cdot |_\Phi $ , then $\frak S_{\Phi^{\ast}}$, the set of compact operators $K$
so that $|K|_{\Phi^{\ast}} = sup \{ |KP|_{\Phi^{\ast}} | P \in \sP \} < \infty $ is its dual with respect to the trace-pairing ( \cite{3} ).

\bigskip
\noindent
{\bf Proposition 4.1}
Assume $ \tau = \tau^{\ast} $, that is $T_j = T^{\ast}_j , 1 \le j \le n $ and let $P, Q \in \sP, PQ = 0$ . Let further
\[
{\Omega = \{X_j = X^{\ast}_j , 1 \le j \le n | i \sum_j [T_j , X_j ] \in (\sB (\sH))_+ +  \mathcal C_1, \sum_j |X_j|_{\Phi^{\ast}} = 1 \}.}
\]
Then we have 
\[
{k_{\Phi} (\tau; P, Q) \ge sup \{ Tr PYP - Tr ((I - P - Q)Y(I - P - Q))_{-} | Y = i \sum_j [T_j , X_j] ,(X_j)_{1 \le j \le n} \in \Omega \}}
\] 
and equality holds if $k_{\Phi} (\tau ; P, Q) > 0$.

\bigskip
\noindent
{\bf {\em Proof.}}
We shall first prove $\ge$ and then assuming $k_{\Phi} (\tau ; P, Q) > 0$ we shall prove $\le$, which will yield the equality stated above. We start with $0 \le B \le I - P - Q , B \in \sR^+_1 $ and $(X_j)_{1 \le j \le n} \in \Omega$ and we will show that
\[
{max_{1 \le j \le n} |[P + B , T_j]|_{\Phi} \ge Tr PYP - Tr ((I - P - Q) Y (I - P - Q))_{-} }
\]
where
\[
{ Y = i \sum_j [T_j , X_j] . }
\]
 
We have
\[
\begin{aligned}
{ Tr ((P + B) Y) = Tr PYP + Tr (B ( I - P - Q) Y ( I - P - Q) )} \\{ = Tr PYP + Tr ( B (( I - P - Q)Y(I - P - Q))_+) - Tr(B ((I - P - Q) Y ( I - P - Q))_{-} )} \\ { \ge Tr PYP -  Tr (B((I - P - Q)Y(I - P - Q))_{-} )  = Tr PYP - Tr (Z^{1/2} B Z^{1/2} )}
\end{aligned} 
\]
where
\[
{ Z = (( I - P - Q) Y ( I - P - Q))_{-} .}
\]
Since $0 \le B \le I $ this gives 
\[
{Tr ((P + B)Y)  \ge Tr PYP  - Tr Z .}
\]
On the other hand
\[
\begin{aligned}
{|Tr (P + B)Y|_{\Phi} = |Tr (i (P+B) \sum_j [T_j , X_j])|_{\Phi} = } \\ { |Tr (\sum_j [P+B, T_j] X_j )|_{\Phi}  \le ( max_{1 \le j \le n} |[P+B , T_j]|_{\Phi} ) \sum_j |X_j|_{\Phi^{\ast}} = max_{1 \le j \le n} |[P + B , T_j] |_{\Phi} .}
\end{aligned}
\]
Hence
\[
{ Tr PYP - Tr (( I - P - Q) Y (I - P - Q) )_{-} \le max_{1 \le j \le n} |[ P + B, T_j] |_{\Phi}}
\]
for all $B \in \sR^+_1 , 0 \le B \le I - P - Q $.

This gives
\[
{Tr PYP - Tr (( I - P - Q) Y (I - P - Q))_{-}  \le k_{\Phi} (\tau ; P , Q) }
\]
since $k_{\Phi} (\tau; P , Q) $ is the $inf$ of the $ max_{1 \le j \le n } |[T_j , P+B]|_{\Phi} $ when $B \in \sR^+_1$ satisfies $0 \le B \le I - P - Q $. This concludes the proof of $\ge$.

\bigskip
Assume now that $k_{\Phi} (\tau; P , Q) > 0$. To prove $\le $ we shall consider the real Banach space $(\frak S^{(0)}_{\Phi , h})^n $ of $n$-tuples $(X_j )_{1 \le j \le n}, X_j = X^{\ast}_j \in \frak S^{(0)}_{\Phi} $ with the norm $ max_{1 \le j \le n} |X_j|_{\Phi}$ and two disjoint convex subsets of this Banach space. The first is the open ball centered at $0$ of radius $k_\Phi (\tau; P, Q)$. The second is:
\[
{ \{ (i [T_j , A] )_{1 \le j \le n}  |  A \in P + (I - P - Q) \sR^{+}_1 (I - P - Q) \} .}
\]
The two convex sets are disjoint and the first is open, so that there is 
\[
{ (X_j)_{1 \le j \le n} \in (( \frak S^{(0)}_{\Phi , h} )^n  )^{\ast}  = ( \frak S_{\Phi^{\ast}, h } )^n }
\]
separating the two and having norm 1. Thus we have
\[
{\sum_j |X_j |_{\Phi^{\ast } } = 1 }
\]
and
\[
{\sum_j Tr (i [T_j , A] X_j) \ge k_{\Phi } (\tau; P, Q)}
\]
for all $A \in P + (I - P - Q) \sR^+_1 (I - P - Q) $ ,
This gives
\[
{\sum_j Tr (i [T_j , P+ B] X_j)  \ge k_{\Phi} (\tau ; P , Q) }
\]
for all $ B \in \sR^+_1 , 0 \le B \le I - P - Q $.
The $ LHS $ equals
\[
{ Tr ((P + B ) i \sum_j [X_j , T_j] ) = Tr PYP + Tr (B (I - P - Q) Y(I - P - Q) ) }
\]
where
\[
{ Y = i \sum_j [X_j , T_j] .}
\]
Remark that $(I - P - Q) Y (I - P - Q) - Y $  is a finite rank operator of rank $ \le 2 Tr (P + Q) $ and its norm is $ \le 4 \sum_j \| T_j \| \| X_j \| \le 4 max_{1 \le j \le n } \| T_j \| $. This implies that
\[
{ inf \{ Tr CY | C\in \sR^+_1 \} \ge const + inf  \{ Tr BY | B \in (I - P - Q) \sR^+_1 (I - P - Q) \} >  - \infty }
\]
so that  $ Y \in ( \sB (\sH ) )_{+}  +  \mathcal C_1 $. Thus we have also proved that
\[
{ inf \{ Tr PYP + Tr (B (I - P - Q) Y ( I - P - Q) | B \in \sR^+_1 \} \ge k_{\Phi} ( \tau ; P , Q ) . }
\]
The above $inf$ is precisely 
\[
{ Tr PYP - Tr ((( I - P - Q) Y (I - P - Q ) )_{-} )  }
\]  
which is the result we wanted to prove.
\qed 

\bigskip
There is an analogue of Proposition 4.1 for unitary operators or equivalently for the corresponding inner automorphisms. The proof being along the same lines is left as an exercise for the reader.

\bigskip
\noindent
{\bf Proposition 4.2}
Let $u$ be a $n$-tuple of unitary operators and let $P, Q \in \sP, PQ = 0 $. Let further
\[
{ \Omega = \{ X_j = X^{\ast}_j \in \frak S_{{ \Phi}^{\ast}}, 1 \le \j \le n | \sum_j (( Ad U_j )(X_j) - X_j ) \in \sB (\sH)_{+} + \mathcal C_1, \sum_j |X_j|_{\Phi^{\ast}} = 1 \} .} 
\]
Then we have:
\[
\begin{aligned}
{ k_{\Phi} (u : P, Q) = k_{\Phi} (Ad u ; P, Q) } \\{ \ge sup \{ Tr PYP - Tr ((I - P - Q) Y ( I - P - Q))_{-} | Y = \sum_j ((Ad U_j)(X_j) - X_j), (X_j)_{1 \le j \le n} \in \Omega  \} } 
\end{aligned} 
\]
and if $k_{\Phi} (u; P, Q) > 0 $ equality holds. 

\section{An example}
\label{sec5}

Let $U$ be the bilateral shift operator on $\sH = \ell^2 ( \Bbb Z) , U e_j = e_{j + 1} $ where $\{ e_j \}_{j \in \Bbb Z} $ is the canonical orthonormal basis. If $ f : \Bbb Z \longrightarrow \Bbb C $ is a bounded function we will denote by $ D (f) $ the diagonal operator in $\sH$ with respect to the canonical basis. In case $ f : \Bbb R \longrightarrow \Bbb C $ we will write $D(f) $ for $D (f| \Bbb Z) . $ Moreover if 
$\omega \subset \Bbb Z , P_{\omega} $ will denote the projection $D( \chi_{\omega} ) $ . In case of a singleton $\{j\} $ we will write $P_j$ instead of $P_{\{j\}} $.

Here we will compute 
\[
{k_p (U; P_M , P_N) }
\]
where $ M , N \subset \Bbb Z $ are two disjoint finite nonempty subsets, By Proposition 3.2 and Proposition 3.4 this would be equivalent to a problem in $\ell^{\infty} ( \Bbb Z) $ , that is a problem on the Cayley graph of $ \Bbb Z $. We will not use this explicitly, though all our computations will be around two operators $ B = D(f)$ and $X = D(g) $ .

Let $ a = inf (M \cup N ) $ and $ b = sup ( M \cup N ) $  and let $ a \le a_1 < b_1 \le a_2 < b_2 \le
  \dots \le a_m < b_m \le b $ be so that the 
$ (a_j , b_j ) $ are the maximal open intervals so that $ (a_j , b_j) \cap (M \cup N ) = \varnothing $ and the endpoints $a_j$ and $b_j$ are in different sets of the partition of $M \cup N $ into $M$ and $N$ .  Let also $h \in \Bbb N $. We define a continuous function $f : \Bbb R \longrightarrow \Bbb R$ as follows. First we require that
\[
{ f|M \equiv 1, f|N \equiv 0 , f| ( - \infty, a - h] \equiv 0 , f| [ b + h, \infty ) \equiv 0 .} 
\] 
Then on each open interval at the endpoints of which $f $ has been defined we extend the definition by linearity on the interval. Thus $f$ will
be piecewise linear with respect to the partition $ - \infty < a-h < a \le a_1 < b_1 \le a_2 < b_2 \le  \dots \le a_m < b_m \le b < b + h < \infty $. The intervals on which this function is not constant are the $[a_j , b_j ] $ and possibly also $[a - h, a]$ and $[b , b + h ] $ depending on whether $ a, b \in M$ or not. Thus if $ B = D(f) $ the list of non-zero singular values of $B - Ad(U)(B) $ consists of $ b_j - a_j $ times the number $(b_j - a_j)^{-1} $ for each interval $[a_j , b_j]$ and each of the intervals $ [a - h, a] , [b , b + h] $ may contribute $h$ times the number $h^{-1} $ depending on whether $ a, b \in M $.  This gives:
\[
{| B - Ad (U)(B) |^p_p = \sum_j (b_j - a_j)^{1 - p} + h^{1 - p} \cdot \sharp ( \{ a, b \} \cap M ) .}
\]
Since 
\[
{ | B - Ad (U) (B) |_p \ge k_p ( Ad (U) ; P_M , P_N ) }
\] 
we get the following upper bound.  If $ p =1$ , we have
 \[
                                                                                                                                                                                                                                  {m + \sharp ( \{ a, b\} \cap M ) \ge k_1 ( Ad (U) ; P_M , P_N) }  
\]                                                                                                                                                                                                                                                                                                                                                                                                                                                                                                                                                                                                                                                                                              while if $1 < p < \infty$ we have 
\[
{ (\sum_j (b_j - a_j)^{ 1 - p} ) ^{1/p} \ge k_p ( Ad (U); P_M , P_N ) }
\]
because $h$ being arbitrary we can take the limit as $ h \rightarrow \infty $.

\bigskip
To get the lower bound using Proposition 4.2 we construct an operator $X = D(g) $. Let $\epsilon (j) = -1 $ if $ a_j \in M $ and $b_j \in N $
and let $ \epsilon (j) = + 1$ if $ a_j \in N $ and $ b_j \in M $ and observe that $ \epsilon (j) =  - \epsilon (j + 1) $. If $ p > 1$ we define 
\[
{ g = c \cdot \sum_j \epsilon (j) (b_j - a_j) ^{-p/q} \cdot \chi_{[a_j , b_j) }  }
\]
where 
\[
{ c = ( \sum_j (b_j - a_j)^{1 - p} ) ^{ -1/q} . }
\]
We have
\[
{ | X |^q_q = c^q \cdot \sum_j (b_j - a_j) \cdot (b_j - a_j) ^{-p} = 1 .} 
\]
If $Y = Ad (U)(X) - X $ we have
\[
{ Y = c \cdot \sum_j \epsilon (j) (b_j - a_j )^{-p/q)} ( P_{b_j} - P_{a_j} ) .}
\] 
From Proposition 4.2 we get the lower bound 
\[
{Tr P_M Y P_M - Tr ((I - P_M - P_N) Y (I - P_M - P_N))_{-} .}
\]
Since $b_j , a_j \in (M \cup N) $ the second term in the lower bound is zero, so we need only compute  $Tr P_M Y P_M $. If $\epsilon (j) = -1$ we have $a_j \in M $ and $ b_j \in N $ , while if $\epsilon (j) = +1 $ we have $a_j \in N $ and $ b_j \in M $ . This gives that $Tr P_M (\epsilon (j) (P_{b_j} - P_{a_j} )) = 1 $ for all indices $j $. It follows that 
\[
{ Tr P_M Y P_M  = c \cdot \sum_j ( b_j - a_j ) ^{-p/q} }
\]
and because $ p/q = 1/p$  we have 
\[
{Tr P_M Y P_M = ( \sum_j (b_j - a_j)^{1 - p} ) ^{1/p} }
\]
where we used $ 1 - 1/q = 1/p $. Thus the lower and the upper bound are equal if $p > 1$.

\bigskip
To obtain the lower bound when $ p = 1$ we shall consider
\[
{ g = - \epsilon (1) \chi_{( - \infty , a)} - \epsilon (m) \chi_{[b , \infty)} + \sum_j \epsilon (j) \chi_{[ a_j , b_j) } }
\]
and $ X = D(g)$ . Then 
\[
{ Y = Ad (U)(X) - X = -\epsilon (1) P_a + \epsilon (m) P_b + \sum_j \epsilon (j) (P_{b_j} - P_{a_j}) . }
\]
Again, the lower bound reduces to computing
\[
{Tr P_M Y P_M = Tr P_M ( - \epsilon (1) P_a + \epsilon (m) P_b ) P_M  + m . }
\]
It is also easy to see that 
\[
{ Tr P_M ( - \epsilon (1) P_a + \epsilon (m) P_b ) P_M = \sharp ( \{a , b \} \cap M) }
\]
so that also in this case { the lower and upper bounds we found for $ k_{\Phi} (Ad (U) ; P_M , P_N) $
are equal.

Summing up and using Proposition 3.4, we have proved the following result.

\bigskip
\noindent
{\bf Proposition 5.1}
With the notation introduced above we have
\[
{k_1 ( U; P_M, P_N) = m + \sharp ( \{ a , b \} \cap M ) }
\]
and if $ 1 < p < \infty $ we have
\[
 { k_p ( U; P_M , P_N ) = ( \sum_j ( b_j - a_j)^{1 - p} )^{1/p} . }
 \]

\section{The singular projection and the regular projection}
\label{sec6}
We adapt and generalize to our semifinite setting the facts in \cite {9} about the largest projection on which the quasicentral modulus vanishes.

The following Lemma is based on an argument we used in the proof of Proposition 3.2.

\bigskip
\noindent
{\bf Lemma 6.1}
Let  $\phi$ be the function in property 3. of the $\sJ$-norm. If $P_1, P_2, Q_1, Q_2 \in \sP,P_1Q_1 = P_2Q_2 = 0,   \| T_j \| \le C, 1 \le j \le n $ , and $ \epsilon > 0$ is 
so that 
\[
{ |P_1 - P_2|_1 < \epsilon , |Q_1 -Q_2 |_1  < \epsilon , }
\]
then we have
\[ 
\begin{aligned}
{|k_{\sJ} (\tau ;P_1 ,Q_1) - k_{\sJ} (\tau ;P_2,Q_2)| \le 4C \phi (6 \epsilon),} \\  {|k_{\sJ} (\alpha;P_1 ,Q_1) - k_{\sJ} (\alpha ;P_2,Q_2)| \le 4\phi (6 \epsilon) .} 
\end{aligned}
\]

\bigskip
\noindent
{\bf {\em Proof.}}
Let  $ X \in \sR^+_1$ and let 
\[
\begin{aligned}
{A = (4C)^{-1} [T_j, P_1 + (I - P_1 -Q_1)X(I - P_1 - Q_1] } \\ { B = (4C)^{-1} [T_j ; P_2 + (I - P_2 - Q_2)X(I - P_2 - Q_2) ] .}
\end{aligned}
\]
Then we have $\| A - B \| \le 1 $ and $ | A - B| \le 6\epsilon $ so that 
\[
{ | A - B |_{\sJ} \le \phi (6 \epsilon) .}
\]
Taking into account the way we defined $k_{\sJ} (\tau ; P, Q) $ using a $max$ over $j$ and then an $inf$ over $X$, this gives
\[
{ |k_{\sJ} ( \tau ; P_1 , Q_1) -  k_{\sJ} (\tau ; P_2, Q_2) | \le 4C \phi (6\epsilon) . }
\]

For automorphisms we use the same argument with $A, B$ defined now to be
\[
{ A = 4^{-1} (\alpha_j ( P_1 + (I - P_1 - Q_1) X (I - P_1 - Q_ 1) ) -  (P_1 + (I - P_1 - Q_1) X (I - P_1 - Q_1))}
\]
and
\[
{ B = 4^{-1} (\alpha_j ( P_2 + (I - P_2 - Q_2) X (I - P_2 - Q_2)) - ( P_2 + (I - P_2 - Q_2) X ( I - P_2 - Q_2)) .}
\]
\qed

\bigskip
\noindent
{\bf Proposition 6.1}
Assume that   
\[
{B_m \in \sR , B_m \ge 0,  w-\lim_{m \rightarrow \infty} B_m = B }
\]
and that
\[
{ \lim_{m \rightarrow \infty} \max_{ 1 \le j \le n} | [ T_j , B_m] |_{\sJ} = 0 .}
\]
Then, if $F$ is the support projection of $B$ ( i.e. $E(B ; (0 , \infty)) $) we have $[ T_j , F] = 0, 1 \le j \le n $ and
\[
{ k_{\sJ} ( \tau ; F ) = 0 .}
 \]
 
 \bigskip
 \noindent
 {\bf {\em Proof. }}                                                                                                                                                                                                                                   
Replacing $B_m$ by $FB_m F$, $M$ by $FMF| F \sH$ , $\tau$ by $\tau | F \sH $ etc., it is easily seen that the proof reduces to the case when $F = I$ that is  $Ker B = 0 $. Recall also that $k_{\sJ} (\tau; I) = k_{\sJ} (\tau) $ . So, we need to prove that $ k_{\sJ} ( \tau; P) = 0$ if $P \in \sP$. 

Remark also that we may assume that
\[
{ s - \lim_{m \rightarrow \infty} B_m = B . }
\]
Indeed, we may pass from the initial $B_m$ 's to a subsequence so that 
\[
{ s - \lim_{m \rightarrow \infty}  m^{-1} ( B_{i_1} + \dots + B_{i_m} ) =  B }
\]
and replace $B_m $ by $ m^{-1} ( B_{i_1} + \dots + B_{i_m} ) $ .

Next, we show how to complete the proof when $ B$ is invertible , that is $ E (B ; [0 , \epsilon) ) = 0 $ for some $\epsilon > 0 $ and then go back to the general situation $Ker B = 0$.
Let $h : \Bbb R \longrightarrow [0, 1] $ be a $ C^{\infty}$-function which is $0$ on $ ( - \infty, 0] $ and $1$ on $ [ \epsilon , \infty) $ . Then we have  
\[
{ s - \lim_{m \rightarrow \infty} h(B_m) = h(B) = ! }
\]
and
\[
{ \lim_{m \rightarrow \infty} \max_{1 \le j \le n} | [ h(B_m) , T_j] |_{\sJ} = 0 .}
\]
Hence, replacing $B_m$ by $h(B_m)$, we may assume that $ B_m \in \sR^+_1$ and $ B_m $ converges strongly to $I$. Let $P \in \sP$ .
Remark then that
\[
{ \lim_{m \rightarrow \infty} | B_m - ( P + ( I - P ) B_m (I - P ) |_1 = 0 .}
\]
This follows from $ |P - PB_m P|_1 \rightarrow 0 $ and $ | (I - P) B_m P |_1 \rightarrow 0 $ which in turn follow from the strong convergences
$ P - PB_mP \rightarrow 0 $ and $ (PB_m(I - P) B_m P )^{1/2} \rightarrow 0 $ in the finite von Neumann algebra $ P \sM P | P |\sH $ endowed with the finite trace which is the restriction of $\rho$.

Since also clearly $ \| B_m - (P + (I - P) B_m (I - P) ) \| \le 2 $ we infer that also 
\[
{ | B_m - ( P + (I - P) B_m (I - P) |_\sJ \rightarrow 0 }
\]
as $ m \rightarrow \infty $. This then gives 
\[
{ | [ T_j , P + (I -P ) B_m (I - P ) ] |_\sJ \rightarrow 0 }
\]
which then finally implies 
\[
{ k_\sJ ( \tau; P ) = 0 .}
\]

Returning to the general case, where only $Ker B = 0 $ is assumed, the result we have obtained thus far is easily seen to give that
\[
{ k_\sJ ( \tau ; E ( B , ( \epsilon , \| B \|) ) = 0 }  
\]
if $ \epsilon > 0$ . 
The proof is then completed by observing that given $P \in \sP$  we can find $ P_k \in \sP , P_k \le E( B ; ( 1/k , \| B \| )) $ so that 
\[
{ \lim_{ k \rightarrow \infty} | P_k - P |_1 = 0 }
\]
and use Lemma 6.1 . We may take $P_k$ to be the left support projection of $E(B ; (1/k , \| B \|))P$ , that is the projection onto the closure of the range of this operator.
\qed

\bigskip

There is an entirely analogous result for automorphisms which we record as the next Proposition, the proof of which is omitted,
being only a slight variation on the preceding proof.

\bigskip
\noindent
{\bf Proposition 6.2} 
Assume that  
\[
{ B_m \in \sR , B_m \ge 0,  w-\lim_{m \rightarrow \infty} B_m = B }
\]
and that
\[
{ \lim_{m \rightarrow \infty} \max_{1 \le j \le n} |[\alpha_j (B_m) - B_m |_\sJ = 0 . }
\]
Then if $F$ is the support projection of $B$ ( that is $ E ( B ; (0 , \infty) ) $ we have  $ \alpha_j (F) = F , 1 \le j \le n $ and
\[
{ k_\sJ (\alpha ; F ) = 0 .}
\]

\bigskip
\noindent
{\bf Corollary 6.1 }
Let $P_1, P_2 \in \sP$.
Then we have
\[
{ k_\sJ ( \tau ; P_1 ) = k_\sJ (\tau ; P_2 ) = 0     \Rightarrow k_\sJ (\tau ; P_1 \vee P_2) = 0 } 
\]
and
\[
{ k_\sJ ( \alpha ; P_1) = k_\sJ (\alpha : P_2 ) =0  \Rightarrow k_\sJ ( \alpha ; P_1 \vee P_2) = 0 }.
\]

\bigskip
\noindent
{\bf {\em Proof.}}
We shall prove only the first assertion, the proof of the second being completely analogous. If $k_\sJ (\tau; P_1) = k_\sJ (\tau; P_2) = 0 $ there exist $A_m, C_m \in \sR^+_1 $ and $A , C \in \sM$ so that
\[
\begin{aligned}
{ A_m P_1 = P_1 , C_m P_2 = P_2 } \\ { w-\lim_{m \rightarrow \infty} A_m = A ,  w-\lim_{m \rightarrow \infty} C_m = C ,} \\ {
\lim_{m \rightarrow \infty} |[A_m , T_j ]|_\sJ = 0 , \lim_{m \rightarrow \infty} |[C_m, T_j]|_\sJ = 0, 1 \le j \le n . } 
\end{aligned}
\]
Since $AP_1 = P_1$, we have that $P_1 \sH$ and $Ker A$ are orthogonal. Similarly $P_2 \sH$ and $Ker C $ are orthogonal. On the other hand $ Ker (A + C) = Ker A \cap Ker C $ because $A$ and $C$ are $\ge 0$. Thus $Ker (A + C)$ is orthogonal to $(P_1 \vee P_2) \sH$ . Applying Proposition 6.1 to the sequence $B_m = A_m + C_m $ we get the desired result.
\qed

\bigskip
\noindent
{\bf Proposition 6.3} 
Given $\tau$ there exists a projection $E^0_\sJ (\tau) \in Proj (\sM) $ so that :
\[
{ P \in Proj(\sM) , k_\sJ (\tau ; P) = 0 \Leftrightarrow P \le E^0_\sJ (\tau) .}
\] 
The projection $ E^0_\sJ (\tau)$ is unique, in particular if $\beta$ is an automorphism of $\sM$ which preserves $\rho$ and $\beta (\tau) = (\tau) $, then we have $\beta (E^0_\sJ (\tau)) = E^0_\sJ (\tau) $. 
Similarly, given $\alpha$ there exists a projection $E^0_\sJ (\alpha) \in Proj (\sM) $ so that:
\[
{ P \in Proj (\sM) , k_\sJ (\alpha; P) = 0 \Leftrightarrow P \le E^0_\sJ ( \alpha) .} 
\]
Moreover the projection $E^0_\sJ (\alpha)$ is unique, in particular if $\beta$ is an automorphism of $\sM$ which preserves $\rho$ and $ \beta  \circ \alpha_j  = \alpha_j \circ \beta $, $1 \le j \le n $ , then we have $\beta (E^0_\sJ (\alpha) ) = E^0_\sJ (\alpha) $ .

\bigskip
\noindent
{\bf {\em Proof.}} 
We will only prove the first half of the statement, the arguments being very similarly for the two cases. Moreover, in view of the definition of $k_\sJ (\tau; P ) $ when $P \in Proj (\sM)) $ it is easily seen that what we must prove , is that 
\[
{ E = \bigvee \{ P \in \sP | k_\sJ (\tau ; P) = 0 \} \Rightarrow k_\sJ (\tau ; E ) = 0 .}
\]
The Hilbert space $\sH$ being separable, there is a sequence $P_i, i \in \Bbb N $ so that $ k_\sJ (\tau ; P_i ) = 0, i \in \Bbb N $ and $ E = \bigvee \{ P_i | i \in \Bbb N \} $. Let $ E_i = P_1 \vee \dots \vee P_i $. By Corollary 6.1 we have $ k_\sJ (\tau ; E_i ) = 0 $. Since $E_i \in \sP$ there is $ B_i \in \sR^+_1 $ so that $B_i E_i = E_i $ and $ | [ T_j , B_i ] |_\sJ  < 1/i , 1 \le j \le n $. We can replace the increasing sequence $E_i $ by a subsequence and assume that the $B_i $'s are weakly convergent to some $B$ . Then $BE = E$ and we can apply Proposition 6.1 to infer that $ k_\sJ ( \tau ; E(B ; (0 , \infty )) =0 $ which implies $ k_\sJ (\tau; E ) = 0 $ since
$ E \le E(B ; ( 0 , \infty )) $ .
\qed

\bigskip

We shall call $E^0_\sJ (\tau), E^0_\sJ (\alpha)$ the $\sJ$-singular projection of $\tau$ and respectively $\alpha$. We shall also use the notation $E_\sJ (\tau) = I - E^0_\sJ (\tau), E_\sJ (\alpha) = I - E^0_\sJ (\alpha) $ and call $ E_\sJ (\tau) , E_\sJ (\alpha) $ the 
$\sJ$-regular projection of $\tau$ and $\alpha$ respectively. 

In \cite{9}, in the case of $\sB (\sH) $ and of a normed ideal given by a norming function $\Phi $ we had called a projection $ P $ which is
$\tau$-invariant $\Phi$-well-behaved if $k_{\Phi} (\tau | P \sH ) = 0$. This is equivalent to $k_{\Phi} (\tau ; P ) = 0 $ and we think  that $\Phi$ - singular, the terminology we introduce here, is perhaps a better term for this. 

\bigskip
\noindent
{\bf Corollary 6.2} 
We have 
\[
{ [ T_j , E^0_\sJ ( \tau) ] = 0 , 1 \le j \le n .}
\]
Similarly, we have
\[
{ \alpha_j ( E^0_\sJ (\alpha )) = E^0_\sJ (\alpha), 1 \le j \le n .}
\]

\bigskip
\noindent
{\bf {\em Proof. } }
Also here we will give only the proof of the first assertion, the proof of the second being along the same lines.

It is clear that it suffices to prove that if $ P \in \sP$ is so that $ k_\sJ (\tau ; P ) = 0 $, then there is $P' \in Proj (\sM)$  so that $P' \ge P, k_\sJ (\tau ; P') = 0, [ P', T_j] = 0, 1 \le j , \le n$. Indeed, since $k_\sJ (\tau; P ) = 0$ , there are $B_m \in \sR^+_1$ so that
$ B_m P = P $ and 
\[
{ \lim_{ m \rightarrow \infty} \max_{ 1 \le j \le n } | [ T_j , B_m ] |_\sJ = 0 .}
\]
Passing to a subsequence, we can assume that
\[
{ w - \lim_{ m \rightarrow \infty} B_m = B }
\]
and we will then have $BP = P $ and $[ B, T_j] = 0 , 1 \le j \le n $ . It follows from Proposition 6.1 that $ P' = E ( B ; ( 0 , \infty)) $ has the desired properties.
\qed

\bigskip
\noindent
{\bf Proposition 6.4} 
If $A_m = A^{\ast}_m \in \sR $ are so that $ \|A_m \| \le C $ for all $m \in \Bbb N $ and
\[
{ \lim_{m \rightarrow \infty } \max_{1 \le j \le n} |[ A_m , T_j]|_\sJ = 0 }
\]
then we have
\[
{ s - \lim_{m \rightarrow \infty} A_m E_\sJ (\tau) = 0. }
\] 
Similarly, if  $A_m = A^{\ast }_m \in \sR $ are so that $\| A_m \| \le C $ for all $ m \in \Bbb N $ and
\[
{ \lim_{m \rightarrow \infty} \max_{1 \le j \le n} | \alpha_j (A_m) - A_m |_\sJ = 0 }
\]
then we have
\[
{ s - \lim_{m \rightarrow \infty} A_m E_\sJ (\alpha) = 0 }
\]

\bigskip
\noindent
{\bf {\em Proof. }} 
We will prove only the first assertion, the proof of the second being along the same lines. 
Let $B_m = E_\sJ (\tau) A^2_m  E_\sJ (\tau) $ so that we will have to prove that 
\[
{w - \lim_{m \rightarrow \infty} B_m = 0. } 
\]
Assuming the contrary and passing to a subsequence of this bounded sequence, we will have
\[
{ w - \lim_{m \rightarrow \infty} B_m = B \neq 0 . } 
\]
Then the $B_m$ and $B$ satisfy the assumptions of Proposition 6.1. . It follows that the projection
$E ( B ; (0 , \infty)) \neq 0 $ is so that $ k_\sJ ( \tau ; E( B; ( 0, \infty ))  = 0 $ and hence by Proposition 6.3, $ E(B ; (0, \infty)) \le E^0_\sJ (\tau) $ while obviously $ E( B; (0 , \infty)) \le E_\sJ ( \tau) $. This contradiction concludes the proof.
\qed

\section{Variants}
\label{sec7}
We briefly discuss here modifications of the definition of the condenser quasicentral modulus quantities along lines, which for the quasicentral modulus we already pointed out in  \cite{8} . We use instead of the max comparable devices in the definitions. This is in preparation of the next section where the variants may have some advantages.

\bigskip 
\noindent

Thus Definition 2.1 is modified as follows: 
\[
{ \tilde{k}_\sJ (\tau; P, Q) = \inf \{ | (\sum_{j=1}^n [T_j , A] ^\ast [T_j , A])^\frac12 |_\sJ \vert A  \in \sR^+_1, AP=P, AQ=0 \}} 
\]
\[
{\tilde{k}_\sJ (\alpha ;P, Q) = \inf \{ | (\sum_{j=1}^n ( \alpha _j (A) - A)^2 )^\frac12 |_\sJ \vert A  \in \sR^+_1, AP=P, AQ=0\}} 
\]
This is then extended also to Definition 2.2 and the further $k_\sJ$ quantities are replaced by $\tilde{k}_\sJ$ quantities. Remark that in essence this amounts to replacing 
\[
{\max_{1 \le j \le n} |X_j|_\sJ }
\]
by
\[
\vert
\begin{pmatrix} 
X_{1} \\
\vdots \\ 
X_{n}
\end{pmatrix} \vert_ \sJ ,
\] 
since
\[
{(\sum_{1 \le j \le n} X_j ^\ast X_j )^\frac12 } 
\]
is the positive operator in the polar decomposition of the column operator matrix. We view here the $ 1 \times n $ matrices with entries in $\sM$ as a subspace in $\sM \otimes \frak M_n $ and use item 4. from the properties of $\vert \cdot \vert _\sJ$  in the preliminaries.

\bigskip
\noindent

Though we will mostly use $\tilde{k}_\sJ$ in this paper,  it is also quite natural to consider a Dirac operator construction. This  produces Dirac-versions $ {k}^D_\sJ (\tau ; P,Q)$ etc. where
\[
{{k}^D_\sJ (\tau ;P, Q) = \inf \{| (\sum_{1 \le j \le n} [T_j , A] \otimes e_j |_\sJ \vert A \in \sR^+_1, AP=P, AQ=0\}}
\]
with $e_1, \dots , e_n$ denoting Clifford matrices.

\bigskip
\noindent

We have
\[
{{k}_\sJ (\tau ; P, Q) \le \tilde{k}_\sJ(\tau ; P, Q) \le n {k}_\sJ (\tau ; P, Q) }.
\]
This implies that a ${k}_\sJ$-condenser quantity is zero or infinity iff the corresponding $\tilde{k}_\sJ$-condenser quantity is zero or
respectively infinity. 

\bigskip
\noindent

It is also easy to see that Lemma 2.1. , Proposition 3.1. , Proposition 3.2. and Proposition 3.3 still hold if ${k}_\sJ$ is replaced by $\tilde{k}_\sJ$.
Note however that it may not be the case that Proposition 3.4. remains valid when we pass to $\tilde{k}_\sJ$.

\section{Noncommutative variational remarks}
\label{section 8}
The quasicentral modulus, in its different versions, is based on quantities for which minimization problems can be formulated:
\[
{I_\tau (X) = \max_{1 \le j \le n} |[T_j , X] |_\sJ },
\]
\[
{I_\alpha (X) = \max_{1 \le j \le n} |\alpha_j (X) - X|_\sJ},
\]
\[
{\tilde{I}_\tau (X) = | (\sum_{1 \le j \le n} [T_j , X]^\ast [T_j , X] ) ^\frac 12 |_\sJ} 
\]
\[
{\tilde{I}_\alpha (X) = | (\sum_{1 \le j \le n} (\alpha_j (X) - X )^\ast (\alpha_j (X) - X) )^\frac 12|_\sJ}
\]
\[
{I^D_\tau (X) = |\sum_{1 \le j \le n} [T_j , X] \otimes e_j |_\sJ}
\]
\[
{I^D_\alpha (X) = |\sum_{1 \le j \le n} (\alpha_j (X) - X ) \otimes e_j |_\sJ }
\]
where $\tau = \tau^\ast $ throughout this section.  If $I(X) \in [0 , \infty ] $ denotes any of the above, remark that it is a differential seminorm with additional properties:
\[
{I(X+Y) \le I(X) + I(Y) }
\]
\[
{I(\lambda X) = |\lambda | I(X)}
\]
\[
{I(XY) \le I(X) \Vert Y \Vert + \Vert X \Vert I(Y) }
\]
\[
{I(X) \le C |X|_\sJ }
\]
\[
{ {w} - \lim_{m\to \infty} X_m = X  \Rightarrow  \liminf_{m\to \infty} I(X_m) \ge I(X) }
\]                                                               
and with the exception of $\tilde{I}_\tau$ and $\tilde{I}_\alpha$ we also have 
\[
{I(X^\ast) = I(X) } .
\]
If $ X = X^\ast $ then $\tilde{I}_\tau$ and ${\tilde{I}_\alpha}$ can also be written as follows:
\[
{\tilde{I}_\tau (X) = |(- \sum_{1 \le j \le n} [T_j , X]^2 )^\frac 12 |_\sJ}
\]
\[
{\tilde{I}_\alpha (X) = | (\sum_{1\le j \le n} (\alpha_j (X) - X)^2)^\frac 12 |_\sJ }.
\]
Actually $X = X^\ast $ is a quite natural condition when we set up variational problems.

\bigskip
\noindent

Euler equations in case $\sJ$ is the $p$-class , $2 \le p < \infty$ can be found for the power-scaled $I^p$ when $I$ does not include
a $\max$ . These equations can be viewed as analogues of the $p$-Laplace equation. More precisely let $X = X^\ast $ be such that:
\[
{\frac d{d\epsilon} I^p (X + \epsilon B) \arrowvert_{\epsilon = 0}  = 0}.
\]
for all $B = B^\ast \in \sR $ , where $I(X) < \infty $.

\bigskip                                                                                                                                                
\noindent                                                                                                                                                       
                                                                                                                                                        
  In case $I = \tilde{I}_\tau$ we have
\[
{\frac d{d\epsilon} \rho ((- \sum_{1 \le j \le n} [X + \epsilon B , T_j]^2)^\frac p2 ) \arrowvert_{\epsilon = 0} }
\]
\[
{= \frac p2  \sum_{1 \le k \le n} \rho(-[B,T_k][X,T_k](- \sum_{1\le j \le n} [X,T_j]^2))^{\frac p2 -1} - [B,T_k](- \sum_{1 \le j \le n} [X,T_j]^2)^{\frac p2 - 1} [X,T_k]))}
\]
\[
{= \frac p2 \sum_{1 \le k \le n} \rho (B (-[T_k, [X,T_k] (- \sum_{1 \le j \le n} [X, T_j]^2)^{\frac p2 - 1})  + B (- [T_k, (- \sum_{1 \le j \le n} [X, T_j]^2)^{\frac p2 - 1} [X, T_k] ))}
\]
which gives
\[
{\tag{*}  \sum_{1 \le k \le n} [T_k, [X,T_k] (- \sum_{1 \le j \le n} [X, T_j]^2)^{\frac p2 - 1} + (- \sum_{1 \le j \le n} [X, T_j]^2)^{\frac p2 - 1} [X, T_k]] = 0}.
\]

\bigskip
\noindent

Similarly in case $I = \tilde{I}_\alpha$ we have
\[
{\frac d{d\epsilon} \rho ((\sum_{1 \le j \le n} (\alpha_j (X + \epsilon B) - (X + \epsilon B) )^2 )^{\frac p2} ) \arrowvert_{\epsilon = 0} }
\]
\[
{ = \frac p2 \sum_{1 \le k \le n} \rho ((\alpha_k (B) - B) (\alpha_k (X) - (X)) (\sum_{1 \le j \le n} (\alpha_j (X) - (X))^2)^{\frac p2 - 1} + } 
\] 
\[
            { +  (\alpha_k (B) - B) (\sum_{1 \le j \le n} (\alpha_j (X) - X)^2))^{\frac p2 - 1} (\alpha_k (X) - (X))) } 
\]
\[
{= \frac p2 \sum_{1 \le k \le n} \rho ((\alpha_k (B) - B) D_k) }
\]
\[
{= \frac p2 \sum_{1 \le k \le n} \rho (B (\alpha ^{-1}_k (D_k) - D_k)) }
\]
where
\[
{D_k = (\alpha_k (X) - X)(\sum_{1 \le j \le n} (\alpha_j (X) - X)^2)^{\frac p2 - 1} + (\sum_{1 \le j \le n} (\alpha_j (X) - X)^2))^{\frac p2 - 1} (\alpha_k (X) - X)) }.
\] 
With this notation we have
\[
{\tag {**}  \sum_{1 \le k \le n} (\alpha^{-1}_k (D_k) - D_k) = 0 }
\]

\bigskip 
\noindent

Similar computations can be carried out in the Dirac case.

\bigskip
\noindent
{\bf Remark 8.1}
It is natural to view solutions of $(*)$ as $\tau -p$-harmonic elements and solutions of $(**)$as  $\alpha -p$ -harmonic elements. A possible technical problem which may appear is that in order not to limit considerations to "bounded p-harmonic" elements it may be necessary to be able to handle  the situation when $X$ is an unbounded operator affiliated with $\sM$.

\bigskip
\noindent
{\bf Remark 8.2}
In the case of automorphisms, if $ I \in \sN \subset \sM$ is a von Neumann subalgebra so that $\rho | \sN$ is semifinite and $\alpha_j (\sN) = \sN, 1 \le j \le n$ let $E$ be the conditional expectation of $\sM$ onto $\sN$ so that $\rho \circ E = \rho$. If $I (X) $ stands for $I_{\alpha} (X), \tilde{I}_{\alpha} (X) $ or $I^D_{\alpha} (X)$ , it is easily seen that
\[
{I(X) \ge I(EX)}.
\]

\bigskip
\noindent

The definitions of the condenser quasicentral moduli $ {k}_\sJ (\tau; P, Q), {k}_\sJ (\alpha; P, Q) , \tilde{k}_\sJ (\tau; P, Q) $, $ \tilde{k}_\sJ (\alpha ; P, Q), {k}^D_\sJ (\tau; P, Q), {k}^D_\sJ (\alpha; P, Q) $ where $P, Q \in \sP, PQ = 0$ suggest corresponding variational problems for the $I(X)$ quantities involving the convex sets
\[
{\sC ^0_{PQ} = \{B \in \sR ^{+}_1 | BP = P , BQ = 0 \}}
\]
\[
{\sC _{PQ} = \{X \in \sM | X = X^\ast , 0 \le X \le I, XP = P, XQ = 0 \}}.
\]
The $\inf$ of $I (B)$ when $B \in \sC ^0_{PQ}$ gives the condenser quasicentral moduli, while $\sC ^0_{PQ}$ is weakly dense in $\sC_{PQ}$, which is a weakly compact convex set.

\bigskip

In view of the weak lower semicontinuity property of $I(X)$ we have that the $\inf$ of $I(X)$ over $\sC_{PQ}$ is attained at some point of $\sC_{PQ}$. Note however, that we only know that 
\[
{\inf \{I(X) | X \in \sC_{PQ}\} \le \inf \{I(X) | X \in \sC^0_{PQ}\} }.
\]

\bigskip

Let $X_m \in \sC^0_{PQ} , m \in \bold N $ be a sequence so that
 \[
 {\lim_{m \rightarrow \infty} I(X_m) = inf \{ I(X) | X \in \sC ^0_{PQ} \} }
 \]
 and which is weakly convergent
 \[
 { {w} - \lim_{m \rightarrow \infty} X_m = X_{\infty}}
  \]
 which can be arranged by passing to a subsequence. We have $ X_{\infty} \in \sC_{PQ} $.

 \bigskip
 
 More can be said when $(\sJ , | \cdot |_\sJ )$ is the $p$ - class, $1 < p < \infty$ because then $(\sJ \otimes \frak M _k, | \cdot |_\sJ)$
 is a uniformly convex Banach space. Assume moreover $I(X)$ is one of $\tilde{I} _\tau (X) , \tilde{I} _\alpha (X), I^D_\tau (X), I^D_\alpha (X)$
 i.e. there is no $\max$ in the definition of $I(X)$. Then $I(X) = |\partial (X)|_\sJ$ , where in each of the four cases $\partial (X)$ is :
 \[
 \begin{pmatrix}
 [X, T_1]  \\
 \vdots\\
 [X, T_n]
 \end{pmatrix}
 \]
 \[
 \begin{pmatrix}
 \alpha_1 (X) - X \\
 \vdots \\
 \alpha_n (X) - X \\
 \end{pmatrix} 
 \]
 \[
 \sum_{1 \le j \le n} [X, T_j] \otimes e_j \\
 \]
 \[
 \sum_{1 \le j \le n} (\alpha_j (X) - X) \otimes e_j  
 \] 
 We have 
 \[
 {I (\frac 12 (X_p + X_q)) \ge \lim_{m \rightarrow \infty} I(X_m) }
 \]
 that is
 \[
 {| \frac 12 (\partial (X_p) + \partial (X_q))|_\sJ  \ge \lim_{m \rightarrow \infty} | \partial (X_m) |_\sJ }
 \]
 In view of the uniform convexity we have that the sequence $\partial (X_m) , m \in \bold N $ is convergent in the norm $ | \cdot |_\sJ $.
 Since $\partial (X_m) $ is weakly convergent to $\partial (X_{\infty} )$, we infer that the limit in the norm $|\cdot |_\sJ $ equals the weak limit $\partial (X_{\infty})$. It follows that
 \[
 { I(X_{\infty}) = \inf \{ I(B) | B \in \sC^0_{PQ} \} }
 \]
 
 \bigskip
 
 Assume $X_m^\prime \in \sC^0_{PQ}, m \in \bold N $ is another sequence which is weakly convergent to $X_{\infty} ^\prime$ and so that
 \[ 
 { \lim_{m \rightarrow \infty} I(X_m^\prime) = \inf \{I(B) | B \in \sC^0_{PQ} \}} 
 \]
 Then 
 \[
 {| \frac 12 (\partial (X_m) + \partial (X_m^\prime)) |_\sJ \ge \inf \{ |\partial (B) |_\sJ | B \in \sC^0_{PQ} \} }
 \] 
 which by uniform convexity implies that
 \[
 {\lim_{m \rightarrow \infty} |\partial (X_m) - \partial (X_m^\prime) |_\sJ = 0 }
 \]
 so that
 \[
 {\partial (X_{\infty} ) = \partial (X_{\infty}^\prime ) }
 \]
 The set $\ker \partial $ is a von Neumann subalgebra of $\sM$. In the case of $\tilde{I}_\tau , I^D_\tau$ it is $(\tau)^\prime \cap \sM $ the relative commutant of $\tau$ in $\sM$, while in the case of $\tilde{I}_\alpha , I^D_\alpha $ it is the fixed point algebra of the $n$-tuple of automorphisms $\alpha$.
 
 \bigskip
 
 Summarizing we have shown the following.
 
 \bigskip
 \noindent
 {\bf Remark 8.3}
 Assume $\sJ$ is the $p$-class , $1 < p < \infty$ , and assume that $I(X)$ is one of the four quantities which do not involve a $\max$. Weak limits of $\sC^0_{PQ} $ - sequences which minimize $I$ over $\sC^0_{PQ}$ form a weakly compact convex subset of $\sC_{PQ}$ on which the value of $I$ is the infimum of $I$ over $\sC^0_{PQ}$. Moreover, modulo the von Neumann subalgebra $\ker \partial$ the elements of this convex set are equal.
  
 \bigskip
 \noindent
 
 The preceding remark leaves open the question about equality of the infimum of $I$ over $\sC^0_{PQ}$ and $\sC_{PQ}$. We can answer this in case ${k}_\sJ (\tau) = 0$ or respectively ${k}_\sJ (\alpha) = 0$. There is not a "no $\max$ " restriction on $I(X)$ for this. In this situation there is a sequence $B_k \in \sR^+_1 , k \in \bold N $ so that $B_k (P + Q) = P + Q , s - lim_{k \rightarrow \infty} B_k = I $ and $|[\tau, B_k]|_\sJ \rightarrow 0$ or respectively $|\alpha_j (B_k) - B_k |_\sJ \rightarrow 0 , 1 \le j \le n $ as $k \rightarrow \infty $ that is $I(B_k) \rightarrow 0$ for our choice of $I$. Let $X \in \sC_{PQ} $ be so that
 \[
 {I(X) = \inf \{ I(Y) | Y \in \sC_{PQ} \} }
 \]
Then $B_k X B_k \in \sC^0_{PQ}$ and 
\[
{s - \lim_{k \rightarrow \infty} B_k X B_k = X }
\]
We have
\[
{I(X) \le I(B_k X B_k) \le I(X) + 2 I(B_k) }
\]
This gives
 \[
 {\lim_{k \rightarrow \infty} I(B_k X B_k) = I(X)}
 \]
 
 \bigskip
 
 Summarizing we have the following.
 
 \bigskip
 \noindent
 {\bf Remark 8.4}
 Assume  that ${k}_\sJ (\tau) = 0$ or respectively that ${k}_\sJ (\alpha) = 0$. Then we have
 \[
 { I(X) = \inf \{I(B) | B \in \sC^0_{PQ} \} = \inf \{I(Y) | Y \in \sC_{PQ} \} }
 \]
 for some $X \in \sC_{PQ}  $.
 
 \bigskip
 \noindent
 
 Assume that $\sJ$ is the $p$-class, $2 \le p <  \infty$ and that $I( \cdot) $ involves no $\max $ . Let $X \in \sC_{PQ} $ be so that
 \[
 { I(X) = \inf \{I(Y) | Y \in \sC_{PQ} \} .}
 \]
 Let $P_1 = E(X; \{1\})$ and $Q_1 = E(X; \{0\})$.
If $B \in \sR^+_1$ and $\epsilon \in [-1, 0] $ we have
\[
{ X_{\epsilon} = X + \epsilon ((P_1 - P) + (X - X^2)^{\frac 12}) B ((P_1 - P) + (X - X^2)^{\frac 12}) \in \sC_{PQ} }
\] 
 and hence
 \[
 { \frac d{d \epsilon } I^p (X_\epsilon) \arrowvert_{\epsilon = 0} \le 0. } 
 \]
 The computations preceding Remark 8.1 with $B$ replaced by

\[                                                                                                                                                                                                                              {((P_1 - P) + (X - X^2)^{\frac 12}) B((P_1 - P) + (X - X^2)^{\frac 12}) } 
\]
give
 \[
 { \rho (\Theta ((P_1 - P) + (X - X^2)^{\frac 12}) B ((P_1 - P) + (X - X^2)^{\frac 12}) \le 0 }
 \]
 if $I$ is $\tilde {I}_{\tau} $ or $\tilde {I}_{\alpha} $ and $\Theta$ denotes the quantity appearing in $(*)$ and respectively $ (**)$.
 Since $B \in \sR^+_1$ is arbitrary, this implies 
 \[
 {((P_1 - P) + (X - X^2)^{\frac 12}) \Theta ((P_1 - P) + (X - X^2)^{\frac 12}) \le 0 .}
 \]
 This in turn is equivalent to 
 \[
 {(I - P - Q_1) \Theta (I - P - Q_1) \le 0 .} 
 \]
 
 Similarly if $\epsilon \in [0 ,1] $ we have 
 \[
 {X_{\epsilon} = X + \epsilon ((Q_1 - Q) + (X - X^2)^{\frac 12}) B ((Q_1 - Q) + (X - X^2)^{\frac 12}) \in \sC_{PQ}  }
 \]
 and proceeding as above we get that
 \[
 {(I - P_1 - Q) \Theta (I - P_1 - Q) \ge 0 .}
 \]

Thus in case $ I = \tilde I_{\tau} $ we have

\[
{ \tag {***}  P_1, Q_1 \in Proj (\sM), P \le P_1, Q \le Q_1, P_1 Q_1 = 0, XP_1 = P_1, XQ_1 = 0, }
\]
\[
{(I - P - Q_1) \Theta (I - P - Q_1) \le 0 }
\]
and
\[
{ (I - P_1 - Q) \Theta (I - P_1 - Q) \ge 0 }
\]
where
\[
{\Theta = \sum_{1 \le k \le n} [T_k , [X, T_k] (- \sum_{1 \le j \le n} [X, T_j]^2)^{\frac p2 -1} + (- \sum_{1 \le j \le n} [X , T_j]^2)^{\frac p2 - 1} [X,T_k]] .}
\]

\bigskip
\noindent

 In case $ I = \tilde I_{\alpha}$ we have
 
 \[
 { \tag {****} P_1, Q_1 \in Proj (\sM), P \le P_1, Q \le Q_1, P_1Q _1 = 0, XP_1 = P_1, XQ_1 = 0 }
 \]
 \[
 {( I - P -Q_1) ( \sum_{1 \le k \le n } (\alpha^{-1}_k (D_k) - D_k)) (I - P - Q_1) \le 0}
 \]
 and
 \[
 {(I - P_1 - Q) (\sum_{1 \le k \le n} (\alpha^{-1}_k (D_k) - D_k)) (I - P_1 - Q) \ge 0} 
 \]
 where
 \[
 {D_k = (\alpha_k (X) - X) (\sum_{1 \le j \le n} (\alpha_j (X) - X)^2 )^{\frac p2 - 1} + (\sum_{1 \le j \le n} (\alpha_j(X) - X)^2)^{\frac p2 - 1} (\alpha_k (X) - X)) .} 
 \]
 
 \bigskip
 Similar computations can be carried out in the Dirac case.
 
 \bigskip
 
 Summarizing we have proved the following.
 
 \bigskip
 \noindent
 {\bf Remark 8.5}
 Assume that $\sJ$ is the $p$ - class , $2 \le p < \infty$ and that $X$ is a minimizer of $I( \cdot)$ in $\sC_{PQ}$. If $I$ is $\tilde {I}_{\tau} $ then $X$ satisfies $(***)$ and if $I$ is $\tilde {I}_{\alpha} $ then $X$ satisfies $(****)$. These conditions are compressions of  noncommutative  $\tau - p$ - Laplace and respectively $\alpha - p$ -Laplace inequalities.

%%%%%%%%%%%%%%%%%%%%%%%%%%%%%%%%%%%%%%%%%%%%%%%%%%%%%%%%%%%%%%%%%%%%%%%%%%%%%%%%

%%%%%%%%%%%%%%%%%%%%%%%%%%%%%%%%%%%%%%%%%%%%%%%%%%%%%%%%%%%%%%%%%%%%%%%%%%%%%%%%

\end{document}